\documentclass[12pt]{article}
\usepackage{amsmath,amssymb,epsfig,subfigure,fancybox}
\usepackage{graphicx}
\usepackage{placeins}
\usepackage{color}
\newcommand{\one}{{\rm 1}\hspace*{-0.035in} {\rm l}}

\newcommand{\la}{\lambda}
\newcommand{\lbar}{\overline}

\newcommand{\dl}{\delta}

\newcommand{\ph}{\varphi}

\def\l{\left|}
\def\r{\right|}

\newcommand{\mz}{{m_0}}

\newcommand{\wdt}{\widetilde}
\newcommand{\e}{\varepsilon}
\newcommand{\rr}{{\Bbb R}}
\newcommand{\M}{{\cal M}}
\newcommand{\cd}{(\cdot)}
\newcommand{\nd}{\noindent}

\newcommand{\sg}{{\operatorname{sgn}}}

\def\para#1{\vskip .4\baselineskip\noindent{\bf #1}}
\def\qed{\strut\hfill $\Box$}

\newtheorem{thm}{Theorem}[section]

\newtheorem{lem}[thm]{Lemma}

\newtheorem{rem}[thm]{Remark}

\newcommand{\thmref}[1]{Theorem~{\rm \ref{#1}}}
\newcommand{\lemref}[1]{Lemma~{\rm \ref{#1}}}

\newcommand{\remref}[1]{Remark~{\rm \ref{#1}}}

\def\al{\alpha}
\def\th{\theta}

\newcommand{\beq}[1]{\begin{equation} \label{#1}}
\newcommand{\eeq}{\end{equation}}
\newcommand{\bed}{\begin{displaymath}}
\newcommand{\eed}{\end{displaymath}}
\newcommand{\bea}{\bed\begin{array}{rl}}
\newcommand{\eea}{\end{array}\eed}
\newcommand{\ad}{&\!\!\disp}
\newcommand{\aad}{&\disp}
\newcommand{\barray}{\begin{array}{ll}}
\newcommand{\earray}{\end{array}}

\newcommand{\aleq}[1]{\begin{align}\begin{aligned} #1
\end{aligned}\end{align}}

\def\({\left(}
\def\){\right)}
\def\disp{\displaystyle}

\begin{document}
\title{Sign-Error Adaptive Filtering
Algorithms for Markovian
Parameters\thanks{This research
 was supported in part by the Army Research Office under grant
W911NF-12-1-0223.}}

\author{Araz Hashemi,\thanks{Department of Mathematics, Wayne State
University, Detroit,
MI 48202, araz.hashemi@wayne.edu.}
\and G. Yin,\thanks{Department of Mathematics, Wayne State
University, Detroit, MI 48202, gyin@math.wayne.edu.}
\and Le Yi Wang\thanks{Department of Electrical and Computer
Engineering, Wayne State University, Detroit, MI 48202,
lywang@wayne.edu.} }

\date{December 20, 2012}

\maketitle

\begin{abstract}
Motivated by reduction of computational complexity,
this work develops sign-error adaptive filtering algorithms for estimating time-varying system parameters.
Different from the previous work on sign-error algorithms,
the parameters are time-varying and their dynamics are modeled by a discrete-time
Markov chain.
A distinctive  feature of the algorithms  is the
multi-time-scale framework for characterizing parameter variations and algorithm updating speeds.
This is realized by considering the stepsize of the estimation algorithms and
a scaling parameter that defines the transition rates of the Markov
jump process. Depending on the relative time scales of these two processes, suitably scaled sequences
of the estimates
are shown to converge to either an
ordinary differential equation, or a set of ordinary differential equations modulated by random switching, or
 a stochastic differential equation,
or stochastic differential equations with random switching. Using weak convergence methods, convergence and rates of convergence of the algorithms
are obtained for all these cases.

\vskip 0.1 true in
\nd{\bf Key Words.}
Sign-error algorithms,
regime-switching models,
stochastic approximation, mean squares errors,
convergence, tracking properties.

\vskip 0.1 true in
\nd{\bf EDICS.} ASP-ANAL

\vskip 0.1 true in
\nd{\bf Brief Title.} Sign-Error  Algorithms
for Markovian Parameters

\end{abstract}

\newpage
\setlength{\baselineskip}{0.28in}
\section{Introduction}
Adaptive filtering algorithms have been studied extensively, thanks to their simple
recursive forms and wide applicability for diversified practical problems arising in
estimation, identification, adaptive control, and signal processing \cite{WidrowS}.

Recent rapid advancement in science and technology has introduced many emerging applications in which adaptive
filtering is of substantial utility, including consensus controls, networked systems,
and wireless communications; see \cite{BenvenisteGR,BenvenisteMP,Chen03,
CY03,Eweda,EwedaM,Guo,HMV95,HP98,KushnerS,KushnerY,Kwong,Ljung,MacchiE,PW97b,Ver79,Yin95}.
One typical scenario of such new domains of applications is that the underlying systems are inherently
time varying and their parameter variations are stochastic \cite{YIK09,YinKC,YinKI}.
One important class of such stochastic systems involves systems whose randomly time-varying parameters can be described by Markov chains.
For example, networked systems include communication channels as part of the system topology. Channel connections, interruptions, data transmission queuing and routing, packet delays and losses, are always random. Markov chain models become a natural choice for such systems. For control strategy adaptation and performance optimization, it is essential to capture time-varying system parameters during their operations, which lead to the problems of identifying Markovian regime-switching systems pursued in this paper.

When data acquisition, signal processing, algorithm implementation are subject to resource limitations, it is highly desirable to reduce data complexity. This is especially important when data shuffling involves communication networks. This understanding has motivated the main theme of this paper by  using sign-error updating schemes, which carry much reduced data complexity, in adaptive filtering algorithms, without detrimental  effects on parameter estimation accuracy and convergence rates.

In our recent work, we developed a sign-regressor algorithm for adaptive filters \cite{YHW11}.
The current paper further develops sign-error adaptive filtering algorithms.
It is well-known that sign algorithms have the advantage of reduced computational complexity.
The sign operator reduces the implementation of
the algorithms to bits in data communications and simple
bit shifts in multiplications. As such,  sign algorithms are highly appealing for practical applications.
The work \cite{Gersho} introduced sign algorithms and has inspired much of the subsequent developments in the field.
On the other hand, employing sign operators in adaptive algorithms has introduced substantial challenges in establishing
 convergence properties and error bounds.

 A distinctive  feature of the algorithms introduced in this paper  is the
multi-time-scale framework for characterizing parameter variations and algorithm updating speeds.
This is realized by considering the stepsize of the estimation algorithms and
a scaling parameter that defines the transition rates of the Markov
jump process. Depending on the relative time scales of these two processes, suitably scaled sequences
of the estimates
are shown to converge to either an
ordinary differential equation, or a set of ordinary differential equations modulated by random switching, or
 a stochastic differential equation,
or stochastic differential equations with random switching. Using weak convergence methods, convergence and rates of convergence of the algorithms
are obtained for all these cases.

The rest of the paper is arranged as follows. Section \ref{sec:for} formulates the problems and introduces the two-time-scale framework.
 The main algorithms are presented in Section \ref{sec:MS}. Mean-squares errors on parameter estimators are derived.
  By taking appropriate continuous-time interpolations, Section \ref{sec:con} establishes
convergence properties of  interpolated sequences of estimates from the adaptive filtering
algorithms. Our analysis is based on weak convergence methods.
The convergence properties are obtained by using martingale averaging
techniques. Section \ref{sec:rate} further investigates the rates of convergence. Suitably
interpolated sequences are shown to converge to
either stochastic differential equations or
randomly-switched stochastic differential equations, depending on relations between  the two time scales. Numerical results by simulation are
presented to demonstrate the performance of our algorithms in Section \ref{sec:num}.

\section{Problem Formulation}\label{sec:for}
Let
\beq{ob} y_n = \ph'_n \al_n + e_n, \ n=0,1, \ldots, \eeq
where $\ph_n \in \rr^r$ is the sequence of regression vectors,
$e_n \in \rr$ is a sequence of zero mean random variables representing
the error or noise,
$\al_n \in \rr^r$ is the time-varying true parameter process,
and $y_n \in \rr$ is the sequence of observation signals at time $n$.

Estimates of  $\al_n$ are denoted by $\th_n$ and are given by the following adaptive filtering algorithm
using a sign operator on the prediction error
\beq{alg} \th_{n+1} = \th_n + \mu \ph_n \sg{(y_n - \ph'_n \th_n)} \eeq
where $\sg(y)$ is defined as
$\sg(y) := 1_{\{ y>0\}} - 1_{\{y <0\}}$ for $y \in \rr^1$.
We impose the following assumptions.

\begin{itemize}

\item[(A1)] $\al_n$ is a discrete-time homogeneous
Markov chain with
state space  \beq{statespace} \M= \{ a_1,\ldots,a_\mz\}, \ a_i
\in\rr^r, i=1,\ldots,\mz, \eeq
and whose transition probability matrix is
given by
\beq{tr} P^\e= I+ \e Q,\eeq
where $\e>0$ is a small parameter, $I$ is the
$\rr^{\mz\times \mz}$ identity matrix, and $Q=(q_{ij})\in \rr^{\mz\times
\mz}$ is an irreducible
generator (i.e., $Q$ satisfies $q_{ij} \ge 0$ for $i\not =j$ and
$\sum^\mz_{j=1}q_{ij}=0$ for each $i=1,\ldots,\mz$)
of a continuous-time Markov chain.
 For simplicity,
 assume that the initial distribution
 of the Markov chain
 $\al_n$ is given by $P(\al_0= a_i)= p_{0,i}$, which is
 independent of $\e$ for each $i=1,\ldots,\mz$, where $p_{0,i} \ge 0$
 and $\sum^\mz_{i=1}p_{0,i}=1$.

\item[(A2)]
The sequence of signals $\{ (\ph_n, e_n) \}$ is uniformly bounded, stationary, and
independent of the parameter process $\{ \al_n \}$.
Let ${\cal F}_n$ be the $\sigma$-algebra generated
by  $\{ (\ph_j,e_j), \al_j: j < n; \al_n\}$, and denote the
conditional expectation with respect to ${\cal F}_n$ by $E_n$.

\item[(A3)] For each $i=1,\ldots, \mz$,
define
\aleq{ \label{g-def}
& g_n := \ph_n \sg(\ph'_n[\al_n - \th_n] +e_n)
\\& g_n(\th, i):= \ph_n \sg(\ph'_n[a_i - \th] +e_n)I_{\{\al_n=a_i\}}
\\& \wdt g_n(\th, i) := E_n g_n(\th, i)
}
For each $n$ and $i$, there is an $A_n^{(i)} \in \rr^{r \times r}$ such
that given $\al_n = a_i$, 
\begin{align} \begin{aligned}\label{ani}
& \wdt g_n(\th, i) = A^{(i)}_n (a_i - \th)I_{\{\al_n=a_i\}} + o(|a_i - \th|I_{\{\al_n=a_i\}})
\\& E A^{(i)}_n = A^{(i)}
\end{aligned} \end{align}

\item[(A4)]
There is a sequence of non-negative real numbers $\{ \phi(k) \}$ with
$\sum_{k} \phi^{1/2}(k) < \infty$ such that for each
$n$ and each $j > n$, and for some $K >0$, 
\beq{amix} |E_n A^{(i)}_j - A^{(i)}| \le K \phi^{1/2}(j-n) \eeq
uniformly in $i=1,\ldots, \mz$.
\end{itemize}

\begin{rem}\label{about-assum-sig}{\rm
Let us take a moment to justify the practicality of the assumptions. The boundedness assumption in {\rm (A2)} is fairly mild. For example, we may use a truncated Gaussian process. In addition, it is possible to accommodate unbounded signals by treating martingale difference sequences (which make the proofs slightly simpler).

In {\rm (A3)}, we consider that while $g_n(\th, i)$ is not smooth
w.r.t. $\theta$, its conditional expectation $\wdt g_n(\th, i)$ can be a smooth function of $\th$. The condition \eqref{ani} indicates that $\wdt g_n (\th,i)$ is locally (near $a_i$) linearizable. For example, this is satisfied if the conditional joint density of $(\ph_n, e_n)$ with respect to $\{\ph_j, e_j, j<n, \ph_n\}$ is differentiable with bounded derivatives; see \cite{CY02} for more discussion.
Finally, {\rm (A4)} is essentially a mixing condition which indicates that the remote past and distant future are asymptotically independent. Hence we may work with correlated signals as long as the correlation decays sufficiently quickly between iterates.
}\end{rem}

\section{Mean Squares Error Bounds}\label{sec:MS}
Denote the sequence of estimation errors by $\wdt{\th}_n : = \al_n - \th_n$.
We proceed to obtain bounds for the mean squares error in terms of the
transition rate of the parameter $\e$ and the adaptation rate of the
algorithm $\mu$.

\begin{thm}\label{MSE}
Assume {\rm (A1)}--{\rm (A4)}. Then there is an $N_\e>0$ such that
for all $n \ge N_\e$,
\beq{eq-msq} E | \wdt \th_n |^2 = E |
      \al_n - \th_n|^2 = O\(\mu +\e + \e^2/\mu\) .\eeq
\end{thm}

\para{Proof.}
Define a
function by $V(x) = (x'x)/2$. Observe that
\beq{th-til}
\wdt\th_{n+1} = \al_{n+1} - \th_{n+1} =  \wdt{\th}_n
  - \mu \ph_n \sg(\ph'_n \wdt{\th}_n + e_n)+( \al_{n+1} -\al_{n} )
\eeq
so
\aleq{ \label{v-0}
 E_n V(\wdt\th_{n+1}) - V(\wdt\th_n)\ad =
 E_n\wdt\th'_n[ (\al_{n+1} - \al_n) - \mu\ph_n\sg(\ph'_n \wdt\th_n
+e_n)]\\
 \aad \
 + E_n | (\al_{n+1} - \al_n) - \mu\ph_n\sg(\ph'_n\wdt\th_n + e_n) |^2
}

By {\rm (A2)}, the Markov chain $\al_n$ is independent of $(\ph_n,
e_n)$ and $I_{\{\al_n =a_i\}}$ is ${\cal F}_n$-measurable.
Since the transition matrix is of the form $P^\e= I+ \e Q$, we obtain
\aleq{ \label{mark}
 E_n(\al_{n+1} - \al_n)\ad = \sum^\mz_{i=1} E (\al_{n+1} - a_{i}
\Big| \al_n = a_i) I_{\{\al_n = a_i \} }
\\ \ad = \sum^\mz_{i=1}
\Big[\sum^\mz_{j=1} a_j( \dl_{ij}+ \e q_{ij} ) - a_i
\Big] I_{\{\al_n =a_i\}}
=O(\e)
}
Similarly,
\aleq{ \label{mark-a}
 E_n |\al_{n+1}-\al_{n}|^2
\ad = \sum^\mz_{j=1}\sum^\mz_{i=1} | a_j -a_i|^2
I_{\{\al_n =a_i\}} P(\al_{n+1}=a_j| \al_n=a_i)
\\ \ad = \sum^\mz_{j=1}\sum^\mz_{i=1} | a_j -a_i|^2
I_{\{\al_n =a_i\}}
(\dl_{ij}+ \e q_{ij})\\ \ad =O(\e)
}

Note that $|\wdt \th_n| = |\wdt \th_n| \cdot 1 \le (|\wdt \th_n|^2
+1)/2$, so
\beq{thes}
O(\e) |\wdt \th_n| \le O(\e) (V(\wdt \th_n)+1).\eeq
Since the signals $\{ (\ph_n, e_n)\}$ are bounded, we have
\aleq{\label{thesa}
\ad E_n|(\al_{n+1} - \al_n) - \mu\ph_n\sg(\ph'_n\wdt\th_n +e_n)|^2
\\ \aad \  = E_n|\al_{n+1}-\al_n|^2
+ O(\mu^2 + \mu\e)[V(\wdt\th_n) + 1)
}
Applying (\ref{thesa}) to (\ref{v-0}), we arrive at
\aleq{\label{v-0a}
\ad E_n V(\wdt\th_{n+1}) - V(\wdt\th_n)=
 -\mu E_n \wdt\th'_n\ph_n\sg(\ph'_n\wdt\th_n +e_n)
+ E_n\wdt\th'_n(\al_{n+1} -\al_n)
\\
\aad \quad + E_n|\al_{n+1} -\al_n|^2
+ O(\mu^2 + \mu\e)[V(\wdt\th_n) + 1]
}
Note also that by {\rm (A3)},
\aleq{ \label{v-0b}
 \mu E_n\wdt\th'_n\ph_n\sg(\ph'_n\wdt\th_n +e_n)\ad
= \mu\sum^\mz_{i=1} E_n\wdt\th'_n\ph_n
\sg(\ph'_n\wdt\th_n+e_n) I_{\{\al_n = a_i\}}
\\
\ad = \mu \sum^\mz_{i=1}E_n\wdt\th'_n A^{(i)}_n
\wdt\th_n I_{\{\al_n = a_i\}}
+\mu o(\wdt\th_n)
\\
\ad =\mu \sum^\mz_{i=1}
\wdt\th'_n [A^{(i)}_n -A^{(i)}]
\wdt\th_n I_{\{\al_n = a_i\}}
+\mu \sum^\mz_{i=1}\wdt\th'_n A^{(i)}
\wdt\th_n I_{\{\al_n = a_i\}}
+\mu o(\wdt\th_n)
}

To treat the first three terms in (\ref{v-0a}), we define the following perturbed
Liapunov functions by
\aleq{ \label{v-pert-def}
& V^\mu_1(\wdt\th, n) := \sum^\infty_{j=n}
\sum^\mz_{i=1} -\mu E_n \wdt\th'
[A^{(i)}_j -A^{(i)}] \wdt\th I_{\{\al_j = a_i\}}
\\& V^\mu_2(\wdt\th, n)  := \sum^\infty_{j=n}
\wdt\th' E_n(\al_{j+1} -\al_j)
\\& V^\mu_3(n) := \sum_{j=n}^\infty E_n(\al_{n+1} - \al_n)'(\al_{j+1} -
\al_j)
}

By virtue of {\rm (A4)}, we have
\aleq{\label{v1-small}
|V^\mu_1(\wdt\th,n)| \le\mu \sum^\mz_{i=1} K
|\wdt\th|^2 \sum^\infty_{j=n} \phi^{1/2}(j-n) \le
O(\mu)[V(\wdt\th) +1]
}
Note also that the irreducibility of $Q$ implies that of $I+ \e Q$
for sufficiently small $\e>0$. Thus
there is an $N_\e$ such that for all $ n \ge N_\e$,
$| (I + \e Q)^k - \one \nu_\e| \le \lambda_c^k$
for some $0<\lambda_c<1$, where $\nu_\e$ denotes the
stationary distribution associated with the transition matrix $I+ \e
Q$.
Note that the difference of the $j+1-n$ and $j-n$ step transition
matrices is given by
\bea \ad (I+ \e Q)^{j+1 -n} -(I+\e Q)^{j-n}
 \\
\aad \ =  [ (I+\e Q)-I] (I+\e Q)^{j-n}\\
\aad\  =
[ (I+\e Q)-I]
[(I+\e Q)^{j-n} - \one \nu_\e]  +  [(I+\e
Q)-I]\one \nu_\e \\
\aad \ = (\e Q)[(I+\e Q)^{j-n} - \one \nu_\e ] .\eea
The last line above follows from the fact   $Q \one =0$, hence
$[(I + \e Q)-I] \one \nu_\e  =0$.
Thus
\beq{asy-trans}
\sum^\infty_{j=n} | I+ \e Q)^{j+1 -n} -(I+\e Q)^{j-n}|
\le O(\e) \sum^\infty_{j=n} \lambda^{j-n}_c =O(\e).
\eeq
The forgoing estimates lead to
$\sum^\infty_{j=n} E_n(\al_{j+1}-\al_j) = O(\e)$ and as a result
\beq{v2-small} |V^\mu_2(\wdt \th,n)|
\le O(\e)(V(\wdt \th)+1).\eeq
and similarly
\beq{v3-small}
|V^\mu_3(n)|
=O(\e),\eeq
so all the perturbations can be made small.

Now, we note that
\aleq{ \label{v1-est00}
& E_n V^\mu_1 (\wdt \th_{n+1},n+1)-
V^\mu_1 (\wdt \th_n,n)
\\& = E_n V^\mu_1(\wdt \th_{n+1},n+1)
- E_n V^\mu_1(\wdt \th_n, n+1)
+ E_n V^\mu_1(\wdt \th_n,n+1)-
V^\mu_n(\wdt \th_n,n) .
}
where
\beq{v1-est0}
 E_n V^\mu_1(\wdt \th_n,n+1) -
 V^\mu_1(\wdt \th_n,n) = \mu \sum^\mz_{i=1}
\wdt \th'_n [A^{(i)}_n -A^{(i)}] \wdt\th_n  I_{\{\al_n = a_i\}}
\eeq
and
\aleq{\label{v1-esta}
\ad E_n V^\mu_1 (\wdt \th_{n+1},n+1)
 - E_n V^\mu_1(\wdt \th_n,n+1)\\
 \aad \ =
  \mu \sum^{\infty}_{j=n+1} \sum^\mz_{i=1}
E_n (\wdt \th_{n+1}- \wdt \th_n)'
[A^{(i)}_j -A^{(i)}]  \wdt \th_{n+1} I_{\{\al_n = a_i\}}
\\
\aad \quad + \mu \sum^{\infty}_{j=n+1} \sum^\mz_{i=1}
E_n \wdt \th_n'
[A^{(i)}_j -A^{(i)}] (\wdt
\th_{n+1}- \wdt \th_n)I_{\{\al_n = a_i\}}.
}
Using (\ref{mark}), we have
\aleq{\label{th-est1}
& E_n |\wdt\th_{n+1} - \wdt\th_n|
\le E_n|\al_{n+1} - \al_n| +
\mu E_n|\ph_n\sg(\ph'_n\wdt\th_n+e_n)|
=O(\e +\mu).
}
Thus, in view of {\rm (A4)}
\aleq{\label{v1-estb}
& \l \mu \sum^{\infty}_{j=n+1} \sum^\mz_{i=1}
E_n \wdt \th_n'  E_{n+1}[A^{(i)}_j -A^{(i)}]
(\wdt\th_{n+1}- \wdt \th_n)I_{\{\al_n = a_i\}} \r
\le O(\mu^2 + \mu\e)[V(\wdt\th_n)+1],
}
and
\aleq{\label{v1-estc}
\l \mu \sum^{\infty}_{j=n+1} \sum^\mz_{i=1}
E_n (\wdt \th_{n+1}- \wdt \th_n)'
E_{n+1} [A^{(i)}_j -A^{(i)}]  \wdt \th_{n+1} I_{\{\al_n = a_i\}} \r
\le O(\mu^2 + \mu\e)[V(\wdt\th_n)+1].
}
Putting together (\ref{v1-est00})--(\ref{v1-estc}), we establish that
\aleq{\label{v1-est-fin}
&E_nV^\mu_1(\wdt\th_{n+1}, n+1) - V^\mu_1(\wdt\th_n, n)
=  \mu \sum^\mz_{i=1}
E_n \wdt \th'_n [A^{(i)}_j -A^{(i)}]
\wdt\th_n  I_{\{\al_n = a_i\}}
+ O(\mu^2 + \mu\e)[V(\wdt\th_n)+1].
}
Likewise, we can obtain
\aleq{\label{v2-est}
E_nV^\mu_2(\wdt\th_{n+1}, n+1) - V^\mu_2(\wdt\th_n, n)
= -E_n\wdt\th'_n(\al_{n+1}-\al_n) + O(\e^2 + \mu^2)
}
and
\aleq{\label{v3-est}
E_nV^\mu_3(n+1) -V^\mu_3(n) = - E_n |\al_{n+1}-\al_n|^2 + O(\e^2)
.}

Now we define $$ W(\wdt\th, n) = V(\wdt\th)
+ V^\mu_1(\wdt\th, n) +
 V^\mu_2(\wdt\th, n) + V^\mu_3( n).$$ Since each
$A^{(i)}$ is a stable matrix there is a $\la >0$
such that $\wdt\th' A^{(i)} \wdt\th \ge \la V(\wdt\th)$
for each $i$. Thus we may take $\la$ such that
$-\mu \sum^\mz_{i=1}\wdt\th' A^{(i)}
\wdt\th I_{\{\al_n = a_i\}}
- \mu O(\wdt\th) \le -\la\mu V(\wdt\th)$.
Using this along with (\ref{v-0}), (\ref{v-0b}),
(\ref{v1-est-fin})--(\ref{v3-est}), and the inequality
$O(\mu\e)= O(\mu^2+\e^2)$, we arrive at
\aleq{ \label{w-est0}
& E_n W(\wdt\th_{n+1}, n+1)-W(\wdt\th_n, n)\\
& \qquad
= -\mu\sum^\mz_{i=1}\wdt\th'_n A^{(i)}
\wdt\th_n I_{\{\al_n = a_i\}}
-\mu O(\wdt\th_n)
+O(\mu^2+\e^2)[V(\wdt\th_n)+1]
\\ &\qquad \le -\la\mu V(\wdt\th_n)+
O(\mu^2+\e^2)[V(\wdt\th_n)+1]\\
&\qquad
\le -\la\mu W(\wdt\th_n, n)+O(\mu^2+\e^2)[W(\wdt\th_n,n)+1]
.}

Choose $\mu$ and $\e$ small enough so that there is a $\la_0>0$
satisfying $\la_0 \le \la$ and
$$-\la \mu + O(\mu^2) + O( \e^2) \le -\la_0 \mu.$$
Then we obtain
$$E_n W(\wdt \th_{n+1},n+1) \le (1-\la_0 \mu) W(\wdt \th_n,n)
+ O(\mu^2 + \e^2)
.$$
Note that there is an $N_\e>0$ such that $(1-\la_0 \mu)^{n} \le
O(\mu)$ for $n\ge N_\e$.
Taking expectation in the iteration for $W(\wdt \th_n, n)$ 
and iterating on the resulting inequality yield
\bea \ad E W(\wdt \th_{n+1},n+1) \le (1-\la_0 \mu)^{n} W(\wdt \th_{0},0)
+ O\(\mu + \e^2 /\mu \).\eea 
Thus
$$ E W(\wdt \th_{n+1},n+1) \le O(\mu + \e^2 /\mu)
.$$
Finally, applying (\ref{v1-small})--(\ref{v3-small}) again, we also
obtain
$$E V(\wdt \th_{n+1}) \le O(\mu +\e+ \e^2 /\mu) .$$
Thus the desired result follows. \qed

\section{Convergence Properties}\label{sec:con}

\subsection{Switching ODE Limit: $\mu=O(\e)$}

We assume the adaptation rate and the transition frequency are of
the same order, that is $\mu = O(\e)$. For simplicity, we take $\mu=\e$.
To study the asymptotic properties of the sequence $\{  \th_n\}$, we
take a continuous-time
interpolation of the process.
Define $$  \th^\mu (t)=   \th_n ,\
 \al^\mu (t)= \al_n,\
\hbox{ for }\ t\in [n\mu, n\mu +\mu).$$
We proceed to prove that $  \th ^\mu\cd$ converges
weakly to a system of randomly  switching ordinary differential equations.

\begin{thm}\label{conv}
Assume {\rm (A1)--(A4)} hold and $\e = \mu$. Then
the process $ (  \th^\mu\cd,\al^\mu\cd)$ converges
weakly to $(  \th\cd,\al\cd)$ such that
 $\al\cd$ is a continuous-time Markov chain generated by
 $Q$ and the limit process $  \th\cd$  satisfies
 the Markov switched ordinary differential equation
\beq{sw-ode} \dot   \th(t)= A^{(\al(t))} (\al (t)-   \th(t)), \
\th (0)=  \th_0.\eeq
\end{thm}

The theorem is established through a series of lemmas.
We begin by using a truncation device to bound the estimates.
Define $S_N=\{  \th \in \rr^r: |   \th|\le N\}$ to be the ball with
radius $N$, and $q^N\cd$ as a truncation function that is
equal to 1 for $ \th \in S_N$, 0 for $  \th \in S_{N+1}$,
and sufficiently smooth between.
Then we modify algorithm \eqref{alg} so that
\beq{alg-n}   \th^N_{n+1} :=  \th^N_n + \mu \ph_n
\sg(y_n-\ph'_n   \th^N_n) q^N(  \th^N_n),
\ n=0,1,\ldots, \eeq
is now a bounded sequence of estimates.
As before, define
$$  \th^{N,\mu}(t) :=  \th^N_n
\ \hbox{ for } \ t\in [\mu n,\mu n +\mu).$$

We shall first show that the sequence
$\{\th^{N,\mu}\cd,\al^\mu\cd\}$
is tight, and thus by Prohorov's theorem
we may extract a convergent
subsequence. We will then show the limit satisfies
a switched differential equation. Lastly, we let the truncation bound
$N$ grow
and show the untruncated sequence given by (\ref{alg}) is also weakly
convergent.

\begin{lem}\label{tight}
The sequence $(  \th^{N,\mu}\cd,\al^\mu\cd)$ is
tight in $D([0,\infty):\rr^r
\times \M)$.\end{lem}

\para{Proof of \lemref{tight}.}
Note that the sequence $ \al^\mu\cd $ is
tight by virtue of  \cite[Theorem 4.3]{YinZ05}.
In addition, $\al^\mu\cd$
converges weakly to a Markov chain generated by $Q$.
To proceed, we examine the
asymptotics of the sequence $ \th^{N, \mu}\cd$.
We have that for any $\delta>0$, and $t,s>0$ satisfying $s\le \dl$,
\aleq{\label{t-tight}
&\!\!\! E^\mu_t \l   \th^{N,\mu}(t+s)-  \th^{N,\mu}(t) \r^2
 \le E^\mu_t \l \mu \sum^{(t+s)/\mu -1}_{k=t/\mu}
\ph_k \sg(y_k-\ph'_k   \th^N_k) q^N(  \th^N_k) \r^2
\\& \le \mu^2 E^\mu_t  \sum^{(t+s)/\mu -1}_{j=t/\mu}
\sum^{(t+s)/\mu -1}_{k=t/\mu}
\ph'_j \ph_k \sg(y_j-\ph'_j \th^N_j) \sg(y_k-\ph'_k\th^N_k)
q^N(\th^N_j) q^N(\th^N_k)
\\& \le \mu^2 \sum^{(t+s)/\mu -1}_{j=t/\mu}
\sum^{(t+s)/\mu -1}_{k=t/\mu}
E^\mu_t \l \ph_j \r^2 E^\mu_t\l \ph_k \r^2
\le O(s^2) \le O(\dl^2).
}
For any $T<\infty$ and any $0\le t\le T$, use $E^\mu_t$ to denote the conditional expectation
w.r.t. the $\sigma$-algebra ${\cal F}^\mu_t$,
  we have $$\lim_{\dl\to 0}\limsup_{\mu\to 0} \Big\{\sup_{0\le s\le \delta}
 E [E^\mu_t \l
   \th^{N,\mu}(t+s)-  \th^{N,\mu}(t)\r^2]\Big\} =0.$$
Applying the criterion
\cite[p.47]{Kushner84}, the tightness is proved. \qed

Since $(   \th^{N,\mu}\cd,\al^\mu\cd)$ is tight,
it is sequentially compact. By virtue of Prohorov's
theorem, we can extract a weakly convergence subsequence.
Select such a subsequence and still denote it by $(
\th^{N,\mu}\cd,\al^\mu\cd)$ for notational simplicity. Denote the limit
by
$(  \th^N\cd,\al\cd)$. We proceed to characterize the limit
process.

\begin{lem}\label{weak-t}
The sequence $(  \th^{N,\mu}\cd,\al^\mu\cd)$ converges weakly to
$(\th^N\cd, \al\cd)$ that is a solution of the
martingale problem with operator
\beq{l1-def}
L_1^N f(  \th^N, a_i)  := \nabla f'(  \th^N, a_i)
A^{(i)} [a_i -   \th^N] q^N(  \th^N) + \sum^\mz_{j=1} q_{ij} f(
\th^N, a_j) ,\eeq
where
for each $i\in \M$, $f(\cdot,i)\in C^1_0$ $(C^1$ functions with compact
support$)$.
\end{lem}

\para{Proof.}
To derive the martingale limit, we need only show that
for the $C^1$ function with compact support $f(\cdot,i)$,
 for each bounded and continuous function $h\cd$,
each $t,s>0$, each positive integer $\kappa$,
and each $t_i \le t$ for $i\le \kappa$,
\aleq{\label{mg-p}
\disp E h( \th^{N} (t_i), \al (t_i): i\le \kappa) &
\Big[ f(  \th^N(t+s),\al(t+s))- f(  \th^N (t), \al(t))
\\ &\quad  - \int^{t+s}_t L_1^N
f(  \th^N(\tau),\al(\tau)) d\tau\Big] =0.
}
To verify \eqref{mg-p}, we use the processes indexed by $\mu$.
As before, note that
\aleq{ \label{expan}
\th^{N,\mu} (t+s) -   \th^{N,\mu}(t)
 = \sum^{(t+s)/\mu-1}_{k=t/\mu}
\mu \ph_k \sg(\ph'_k[\al_k -\th_k]+e_k)q^N(\th^N_k).
}
Subdivide the interval
with the end points $t/\mu$ and $(t+s)/\mu-1$ by choosing $m_\mu$ such
that $m_\mu\to \infty$ as $\mu\to 0$ but $\dl_\mu= \mu m_\mu\to 0$.
By the smoothness of $f(\cdot,i)$, it is readily seen that as $\mu\to 0$,
\aleq{ \label{w-es1}
&\!\!\! E h( \th^{N,\mu} (t_i), \al^\mu (t_i): i\le \kappa)
\left[ f(  \th^{N,\mu}(t+s),\al^\mu(t+s))-
f(  \th^{N,\mu} (t), \al^\mu(t)) \right]
\\
& \to Eh( \th^{N} (t_i),
\al (t_i): i\le \kappa)
\left[ f(  \th^N(t+s),\al(t+s))-
f(  \th^N (t), \al(t)) \right].
}
Next, we insert a term to examine the change in the
parameter $\al$ and the estimate $\th^N$ separately
\aleq{ \label{w-es2}
&\!\!\! \lim_{\mu\to 0}
E h( \th^{N,\mu} (t_i), \al^\mu (t_i): i\le \kappa)
\left[ f(  \th^{N,\mu}(t+s),\al^\mu(t+s))-
f(  \th^{N,\mu} (t), \al^\mu(t)) \right]
\\
&\ =\lim_{\mu\to 0}
E h( \th^{N,\mu} (t_i), \al^\mu (t_i): i\le \kappa)
\left[
\sum^{t+s}_{l\dl_\mu=t}[ f( \th^{N}_{lm_\mu+m_\mu},
\al_{lm_\mu+m_\mu})-f( \th^{N}_{lm_\mu},
\al_{lm_\mu})]\right]
\\
&\  =\lim_{\mu\to 0}
E h( \th^{N,\mu} (t_i), \al^\mu (t_i): i\le \kappa) \Big[
\sum^{t+s}_{l\dl_\mu=t}[ f( \th^{N}_{lm_\mu+m_\mu},
\al_{lm_\mu+m_\mu})-f( \th^{N}_{lm_\mu+m_\mu},
\al_{lm_\mu})]
\\
&\hspace*{1.2in} +
\sum^{t+s}_{l\dl_\mu=t}[ f( \th^{N}_{lm_\mu+m_\mu},
\al_{lm_\mu})-f( \th^{N}_{lm_\mu},
\al_{lm_\mu})] \Big].
}

First, we work with the last term in (\ref{w-es2}).
By using a Taylor expansion on each interval indexed by $l$ we have
\begin{align} \begin{aligned}\label{ode0}
&\!\!\! \lim_{\mu\to 0} E h( \th^{N,\mu} (t_i), \al^\mu (t_i): i\le \kappa)
\Big[
\sum^{t+s}_{l\dl_\mu=t}[ f( \th^{N}_{lm_\mu+m_\mu},
\al_{lm_\mu})-f( \th^{N}_{lm_\mu}, \al_{lm_\mu})] \Big]
\\
& = \lim_{\mu\to 0} E h( \th^{N,\mu} (t_i), \al^\mu (t_i): i\le \kappa)
\sum^{t+s}_{l\dl_\mu=t} \Big[
 \dl_\mu {1\over m_\mu} \sum^{lm_\mu+m_\mu-1}_{k=lm_\mu}
\nabla f'( \th^N_{lm_\mu}, \al_{lm_\mu})
 \\
 &  \qquad\quad \hfill  \times \ph_k\sg(\ph'_k( \al_k -\th^N_k) + e_k)q^N(\th^N_k)
\\
& \qquad \qquad \qquad \qquad
+\sum^{lm_\mu+\mu-1}_{k=lm_\mu}
[ \nabla f'(\th^{N,+}_{lm_\mu},\al_{lm\mu}) -
\nabla f'(\th^{N}_{lm_\mu},\al_{lm\mu}) ]
(  \th^N_{k+1}-  \th^N_k)
q^N( \th^N_k)
 \Big].
\end{aligned} \end{align}
where $\th^{N,+}_{lm_\mu}$ is a point on the line segment joining
$\th^N_{lm_\mu}$ and $\th^N_{lm_\mu + m_\mu}$.
Since $$ | \th^N_{lm_\mu+m_\mu} - \th^N_{lm_\mu} | = O(\dl_\mu)$$ and
$\nabla f'(\cdot, i)$ is smooth, we have the last term in \eqref{ode0}
is $o(1)$ in the sense of in probability as $\mu \to 0$. To work with the first term we insert the
conditional expectation $E_{k}$ and apply (\ref{ani}) to obtain
\begin{align} \begin{aligned}\label{ode1}
& \lim_{\mu\to 0} E h( \th^{N,\mu} (t_i), \al^\mu (t_i): i\le \kappa)
 \sum^{t+s}_{l\dl_\mu=t}  \dl_\mu {1\over m_\mu}
 \sum^{lm_\mu+m_\mu-1}_{k=lm_\mu}
 \nabla f'( \th^N_{lm_\mu}, \al_{lm_\mu})\\
 & \qquad \qquad \times
E_{k} [\ph'_k \sg(\ph'_k( \al_k -\th^N_k) + e_k)]q^N(\th^N_k)
\\ & = \lim_{\mu\to 0} E h( \th^{N,\mu} (t_i), \al^\mu (t_i): i\le \kappa)
 \sum^{t+s}_{l\dl_\mu=t}  \sum^{m_0}_{j=1} \dl_\mu {1\over m_\mu}
\sum^{lm_\mu+m_\mu-1}_{k=lm_\mu}
 \nabla f'( \th^N_{lm_\mu}, \al_{lm_\mu})\\
 &\qquad\qquad \times \Big[ A_k^{(j)}(a_j-\th^N_k)
+ o(|\al_k - \th^N_k|) \Big]
 q^N(\th^N_k)I_{\{ \al_k=a_j \} }.
\end{aligned} \end{align}
Then  for small $\mu$,
\beq{ode-es1} E {1 \over m_\mu} \sum^{lm_\mu+m_\mu-1}_{k=lm_\mu}
 \nabla f'( \th^N_{lm_\mu}, \al_{lm_\mu}) o(|\al_k - \th^N_k|)
\le K E|\al_{lm_\mu} - \th_{lm_\mu}| =O(\mu^{1/2}).
\eeq
Letting $\mu l m_{\mu} \to \tau$, then by (\ref{amix}),
\begin{align} \begin{aligned}\label{ode2}
&\!\!\! \lim_{\mu\to 0} E h( \th^{N,\mu} (t_i), \al^\mu (t_i): i\le \kappa)
 \sum^{t+s}_{l\dl_\mu=t}  \sum^{m_0}_{j=1}
 {\dl_\mu \over m_\mu} \sum^{lm_\mu+m_\mu-1}_{k=lm_\mu}
 \nabla f'( \th^N_{lm_\mu}, \al_{lm_\mu})
 A_k^{(j)}(a_j-\th^N_k)\\
 & \qquad \qquad \times q^N(\th^N_k)I_{\{ \al_k=a_j \} }
\\
& = \lim_{\mu\to 0} E h( \th^{N,\mu} (t_i), \al^\mu (t_i): i\le \kappa)
 \sum^{t+s}_{l\dl_\mu=t}  \sum^{m_0}_{j=1}
 {\dl_\mu \over m_\mu} \sum^{lm_\mu+m_\mu-1}_{k=lm_\mu}
 \nabla f'( \th^N_{lm_\mu}, \al_{lm_\mu})
 \Big[ A^{(j)}(a_j-\th^N_k)
 \\
 & \hspace*{3.5in} + [A_k^{(j)} - A^{(j)}](a_j-\th^N_k)
 \Big] q^N(\th^N_k)I_{\{ \al_k=a_j \} }
\\
&=  \lim_{\mu\to 0} E h( \th^{N,\mu}
 (t_i), \al^\mu (t_i): i\le \kappa)
 \sum^{t+s}_{l\dl_\mu=t}  \sum^{m_0}_{j=1}
 \dl_\mu \nabla f'(\th^N_{lm_\mu},\al_{lm_\mu})
 {1 \over m_\mu} \sum^{lm_\mu+m_\mu-1}_{k=lm_\mu}
 A^{(j)} (a_j - \th^N_k)\\
 & \qquad \qquad \times q^N(\th^N_k)I_{\{ \al_k=a_j \} }
\\
& =E h( \th^{N} (t_i), \al (t_i): i\le \kappa)
 \int^{t+s}_t \nabla f' (  \th^N(\tau),\al(\tau))
 A^{(\al(\tau))}[\al(\tau)-\th^N(\tau)] q^N(\tau) d\tau.
\end{aligned} \end{align}
Likewise, we can obtain
\beq{w-es4}\barray
\ad
\lim_{\mu\to 0}
E h( \th^{N,\mu} (t_i), \al^\mu (t_i): i\le \kappa)
\sum^{t+s}_{l\dl_\mu=t}[ f( \th^{N}_{lm_\mu+m_\mu},
\al_{lm_\mu+m_\mu})-f( \th^{N}_{lm_\mu+m_\mu},
\al_{lm_\mu})]\\
\aad \ =
E h( \th^{N,\mu} (t_i), \al^\mu (t_i): i\le \kappa)
\sum^{lm_\mu+m_\mu-1}_{k=lm_\mu}[ f(  \al^N_{lm_\mu},
\al_{lm_\mu+m_\mu})- f(  \th^N_{lm_\mu},\al_{lm_\mu})]
\\
\aad \ =
E h( \th^{N} (t_i), \al (t_i): i\le \kappa)
\left[ \int^{t+s}_t Q f(  \th^N(\tau),\al(\tau)) d\tau\right].\earray\eeq
Combining (\ref{ode0})--(\ref{w-es4}) with (\ref{w-es2}), we have
established \eqref{mg-p} as desired, completing the proof of the lemma.
\qed

\para{Completion of the Proof of \thmref{conv}.}
From \lemref{weak-t}, we have the truncated sequence
$\th^N\cd$ satisfies the switched ODE
$\dot \th^N (t) = A^{(\al(t))}[\al(t) - \th^N(t) q^N(t)], \
\th (0)=  \th_0. $
Next, letting $N\to \infty$,
we show that the
limit of the untruncated
sequence $\th\cd$
and the limit of $\th^N \cd$
as $N\to \infty$ are the same.
The argument is similar to that of
\cite[pp. 249-250]{KushnerY}; we explain the main steps below.
Let $P^0\cd$ and $P^N\cd$ be the
measures induced by $  \th\cd$ and $  \th^N\cd$,
respectively.
Since the martingale problem with operator $L_1^N$ has
a unique solution, the associated differential equation has a
unique solution for each initial condition and
$P^0\cd$ is unique. For each $T<\infty$ and $t\le T$,
$P^0\cd$ agrees with $P^N\cd$ on all Borel subsets of the set of
paths in $D[0,\infty)$ with values in $S_N$.
By using
$P^0(\sup_{t\le T} |  \th (t)|\le N)\to 1$ as $N\to
\infty$,
and the weak convergence of $  \th^{N,\mu}
\cd$ to $  \th^N\cd$,
we have $ \th^\mu\cd$  converges weakly to $  \th\cd$.
Thus the proof of \thmref{conv} is completed. \qed

\begin{rem}\label{f-s-m}
{\rm  The following calculation will be used for both the slow and
fast Markov chain cases. The result is essentially one about
two-time-scale Markov chains considered in \cite{YinZ05}.
Define a probability vector by
$p^\e_n = ( P(\al_n= a_1), \ldots, P(\al_n =a_\mz))\in
\rr^{1\times \mz}.$
Note that $p^\e_0=(p_{0,1}, \ldots, p_{0,\mz})$ (independent of $\e$).
Because the Markov chain is time homogeneous, $ (P^\e)^n$ is the
$n$-step transition probability matrix with
$P^\e=I+\e Q$.
Then,
for  some $0<\la_1<1$,
\beq{appro} p^\e_n = p(\e n) + O(\e + \la_1^{ - n}), \ \  0
\le n\le
O( 1/\e), \eeq
where
$ p(t)=
(p_1(t),\ldots, p_\mz(t)) $ is the
 probability vector  of the
continuous Markov chain with generator $Q$ such that
for all $t \ge 0$
\beq{diff-pi}
 { d p(t) \over dt} = p(t) Q, \\
 \quad p(0) =p_0, \eeq
 and $p_0$ is the initial probability.
 In addition,
 \beq{appro-m} (P^\e)^{n-n_0} =\Xi(\e n_0, \e n) + O(\e + \la_1^{
 - (n-n_0)}),
 \eeq
 where with $t_0 = \e n_0$ and $ t = \e n$, $\Xi(t_0,t)$
 satisfies
 \beq{diff-xi} \left\{ \barray \ad { d \Xi(t_0,t) \over dt}
 =\Xi(t_0,t) Q,\\
 \ad \Xi(t_0,t_0)=I.\earray\right. \eeq
Define the
continuous-time interpolation $\al^\e (t)$
of $\al_n$ as
\beq{eq:thetacont}
\al^\e (t) := \al_n  \ \hbox{ for } \ t \in [n \e ,n
\e  +\e ).    \eeq
 Then $\al^\e \cd$ converges weakly to $\al\cd$, which is
 a continuous-time Markov chain generated by $Q$ with state space
 $\M$.
 The $E \al_n$ can be approximated by
\bea \ad E \al_n = \lbar \al_*(\e n)+
 O(\e +\la_1^{-  n}), \
 \hbox{ for } \ n\le O(1/\e) , \\
\ad \lbar \al_*(\e n)
:= \sum^\mz_{j=1} a_j p_j( \e n) .\eea
}\end{rem}

\subsection{Slowly-Varying Markov Chain: $\e \ll \mu$}
In this case, since the Markov chain
changes so slowly, the time-varying parameter
process is essentially a constant.
 To facilitate the discussion and to
fix notation, we take $\e = \mu^{1+\Delta}$
for some $\Delta>0$  in what follows.

The analysis is similar to the $\e = O(\mu)$ case. Begin
by defining the continuous time interpolation as before. While a truncation device is
still needed, we omit it and assume the iterates are bounded for notational brevity.
The tightness of $\{\th^\mu\cd\}$ can be verified similar to  Lemma \ref{tight}. To
characterize the weak limit we note that the estimates from the previous section
remain valid, except that involving the Markov chain $\al_k$.
Thus we need only examine (from the second to last line of (\ref{ode2}))
\aleq{\label{s-m-1}
& \lim_{\mu\to 0} E h( \th^{N,\mu}
 (t_i), \al^\mu (t_i): i\le \kappa)
 \sum^{t+s}_{l\dl_\mu=t}  \sum^{m_0}_{j=1}
 \dl_\mu \nabla f'(\th^N_{lm_\mu},\al_{lm_\mu})
 \sum^{lm_\mu+m_\mu-1}_{k=lm_\mu}
 {1 \over m_\mu}A^{(j)} a_j I_{\{ \al_k=a_j \} }
\\
&  =\lim_{\mu\to 0} E h( \th^{N,\mu}
 (t_i), \al^\mu (t_i): i\le \kappa)
 \sum^{t+s}_{l\dl_\mu=t}  \sum^{m_0}_{j=1}
 \dl_\mu \nabla f'(\th^N_{lm_\mu},\al_{lm_\mu})
 {1 \over m_\mu}  \sum^{lm_\mu+m_\mu-1}_{k=lm_\mu}
A^{(j)} a_j E_{lm_\mu}I_{\{ \al_k=a_j \} }
 \\
 &  =\lim_{\mu\to 0} E h( \th^{N,\mu}
 (t_i), \al^\mu (t_i): i\le \kappa)
 \sum^{t+s}_{l\dl_\mu=t}  \sum^{m_0}_{j=1}
 \dl_\mu \nabla f'(\th^N_{lm_\mu},\al_{lm_\mu})
 \\
 & \qquad  \hfill \times
{1\over m_\mu}\sum^{lm_\mu+m_\mu-1}_{k=lm_\mu}
\sum^\mz_{i_1=1} \sum^\mz_{i_0=1}
 A^{(j)}  a_{j} P(\al_k =a_{j}| \al_{lm_\mu}=a_{i_1}) P(\al_{lm_\mu}
=a_{i_1}| \al_0=a_{i_0})
P(\al_0=a_{i_0})
\\
& \to E h( \th^{N} (t_i), \al (t_i): i\le \kappa)
\sum^\mz_{i_0=1} \int^{t+s}_t \nabla f' (  \th(\tau),\al(\tau)) A^{(i_0)} a_{i_0}
P(\al_0= a_{i_0}) d\tau.
 }
To obtain the last line above,
we have used that for $ lm_\mu \le k\le l m_\mu + m_\mu$ since  $\e =\mu^{1+\Delta}$,
$\e l m_\mu + m_\mu \le \mu^\Delta(t+s) + \dl_\mu \to 0$ as $\mu \to 0$,
we have by \remref{f-s-m},
\bea  (P^\e)^{k-lm_\mu} \ad = \Xi ( \e k, \e lm_\mu )+ O\(\e +\la_1^{ -
( k-lm_\mu)}\)\\
\ad \to I \ \hbox{ as } \ \mu\to 0 ,\\
(P^\e)^{lm_\mu} \ad = \Xi ( 0, \e lm_\mu )+ O\(\e +
\la_1^{- lm_\mu} )\)\\
\ad \to I \ \hbox{ as } \ \mu\to 0.\eea
We omit the details, but present the main result as follows,

\begin{thm} \label{slow-ch}
Assume {\rm(A1)--(A4)} hold, and $\e =\mu^{1+\Delta}$ for some $\Delta >0$.
Then we have $ \th^\mu\cd$ converges weakly to $  \th\cd$ such that
$ \th\cd$ is the unique solution of the differential equation
\beq{ode-sl} {d\over dt}   \th(t)=  \sum^\mz_{i=1} A^{(i)} (a_i - \th(t))
P(\al_0 = a_i), \   \th(0)=  \th_0. \eeq
\end{thm}

\subsection{Fast-Varying Markov Chain:  $\mu \ll \e$}
The idea for the fast varying chain is that the parameter changes so fast
that it quickly approaches the
stationary distribution of the Markov chain.  As a result, the limit
dynamic system is one that is averaged
out with respect to the stationary distribution of the Markov chain.
In this section, we take $\e = \mu^{\gamma}$ where $ 1/2 < \gamma <1$.
Then, letting $\mu lm_\mu \to \tau$ as in the proof of \thmref{conv}, we
have
$\e (k- lm_\mu)= \mu^\gamma (k-lm_\mu ) \to \infty$. Thus,
for some $0<\la_1<1$,
\bea \ad
\Xi_{ij} (\e lm_\mu,\e k)=  \nu_j + O( \e  + \la_1^{- ( k-
lm_\mu)}),\eea
where
$\nu=(\nu_1,\dots,\nu_{m_0})$ is the stationary distribution of the
continuous-time Markov chain
with generator $Q$,
$\Xi_{ij}(s_1, s_2)$ denotes the
$ij$th entry of the matrix
$\Xi (s_1,s_2)$.
Therefore, we can show that as $\mu\to 0$,
\aleq{\label{f-e-es1}
& \lim_{\mu\to 0} E h( \th^{N,\mu}
 (t_i), \al^\mu (t_i): i\le \kappa)
 \sum^{t+s}_{l\dl_\mu=t}  \sum^{m_0}_{j=1}
 \dl_\mu \nabla f'(\th^N_{lm_\mu},\al_{lm_\mu})
 {1 \over m_\mu}  \sum^{lm_\mu+m_\mu-1}_{k=lm_\mu}
A^{(j)} a_j E_{lm_\mu}I_{\{ \al_k=a_j \} }
\\
&  \to  E h( \th^{N}
 (t_i), \al(t_i): i\le \kappa)
\sum^\mz_{j=1}  \int_t^{t+s}
  \nabla f' (  \th(\tau),\al(\tau)) A^{(j)} a_{j} \nu_j d\tau.
}

\begin{thm} \label{fast-ch}
Assume \rm{(A1)--(A4)} hold, and $\e =\mu^\gamma$ for some
$1/2<\gamma <1$.
Then we have $ \th^\mu\cd$ converges
weakly to $  \th\cd$ such that $ \th\cd$ is the unique solution of
the differential equation
\beq{ode-sl-a} {d\over dt}   \th(t)= \sum^\mz_{j=1} A^{(j)}
( \nu_i a_i -\th(t)), \   \th(0)=  \th_0. \eeq
\end{thm}

\section{Rates of Convergence}\label{sec:rate}
\subsection{Scaled Errors: $\e=\mu$}

Define $u_n := \wdt \th_n / \sqrt{\mu} = (\al_n - \th_n)/ \sqrt \mu$. Then
\beq{uerr-dif}
u_{n+1} = u_n -
\sqrt{\mu} \ph_n \sg(\ph'_n \wdt \th_n +e_n)
+{\al_{n+1} -\al_n \over \sqrt \mu}
\eeq

In view of \thmref{MSE} there is a $N_\mu$ such that
$E | \al_n - \th_n |^2 = O(\mu)$ for $n \ge N_\mu$, with which we can show
$\{ u_n : n\ge N_\mu \} $ is tight.
In addition, take $N_\mu$ large such that by (\ref{asy-trans}), we have
\beq{ord-al}
\sum^\infty_{j=n} E_n(\al_{j+1}-\al_j) = O(\mu) \eeq
Then  define
$$u^\mu (t) := u_n \ \hbox{ for } \ t\in [(n-N_\mu) \mu, (n-
N_\mu)\mu+\mu).$$
We can then  proceed to the study of the asymptotic distribution of
$u^\mu\cd$. As before, a truncation device may be employed.
 For notational simplicity, it will be assumed here.

\begin{lem}\label{u-tight}
The sequence $ \{ u^\mu \cd \} $ is tight in $D( [0,\infty) ; \rr^r)$.
\end{lem}

\para{Proof.}
Note that
\aleq{\label{udiff-1}
& u^\mu(t+s) - u^\mu(t) =
-\sqrt \mu \sum^{(t+s)/\mu -1}_{k=t/\mu}g_k
+ {\al_{(t+s)/\mu} -\al_{t/\mu} \over \sqrt \mu}
.}
Note that we have used the convention that $t/\mu$ denotes the integer part of $t/\mu$
in the above.
Use $E^\mu_t$ to denote the conditional expectation with respect to the
$\sigma$-algebra ${\cal F }^\mu_t = \sigma \{ u^\mu( \tau ) : \tau \le t
\}$.
Then by (\ref{ord-al}),
\aleq{\label{u-tight-00}
& E_t^\mu \l u^\mu(t+s) - u^\mu (t) \r^2
\le K E^\mu_t\l \sum^{(t+s)/\mu -1}_{k=t/\mu} -\sqrt\mu g_k \r^2 + O(\sqrt\mu)
}
Now we examine
\aleq{ \label{u-tight-0}
&
\disp E^\mu_t\l \sum^{(t+s)/\mu -1}_{k=t/\mu} -\sqrt\mu g_k \r^2
 =\sum_{i=1}^\mz E^\mu_t  \mu \sum^{(t+s)/\mu -1}_{k=t/\mu}
 \sum^{(t+s)/\mu -1}_{j=t/\mu}  g_k'  g_j I_{\{\al_k = a_i = \al_j \} }
\\&
 \le \sum^\mz_{i=1}E^\mu_t \mu  \sum^{(t+s)/\mu -1}_{k=t/\mu}
 \sum^{(t+s)/\mu -1}_{j=t/\mu}
[ A_k^{(i)}\wdt\th_k + o(\wdt\th_k)]' [ A_j^{(i)}\wdt\th_j + o(\wdt\th_j)]  I_{\{\al_k = a_i\}}
I_{\{\al_j  = a_i\}}
\\&
\le \sum^\mz_{i=1}E^\mu_t K \mu \l \sum^{(t+s)/\mu -1}_{k=t/\mu}
(A_k^{(i)} - A^{(i)}) \wdt\th_k + o(\wdt\th_k) \r^2 I_{\{\al_k = a_i\}}\\
& \quad
+  \sum^\mz_{i=1}E^\mu_t K \mu \l \sum^{(t+s)/\mu -1}_{k=t/\mu}
A^{(i)}\wdt\th_k \r^2 I_{\{\al_k = a_i\}}
.}
 Since $E|\wdt\th_k|^2 = O(\mu)$ for $k$ large ($\mu$ small), in the last
term of (\ref{u-tight-0}) we have
\aleq{\label{u-tight-1}
E  \sum^\mz_{i=1}E^\mu_t K \mu \l \sum^{(t+s)/\mu -1}_{k=t/\mu}
A^{(i)}\wdt\th_k \r^2 I_{\{\al_k = a_i\}}
\le K \mu \sum^\mz_{i=1}E \sum^{(t+s)/\mu -1}_{k=t/\mu}
\l \wdt\th_k \r^2I_{\{\al_k = a_i\}}
\le O(\mu)s
.}
For the first term we use the mixing inequality of {\rm (A4)},
\aleq{ \label{u-tight-2}
&\!\!\! E \sum^\mz_{i=1}E^\mu_t K \mu \l \sum^{(t+s)/\mu -1}_{k=t/\mu}
(A_k^{(i)} - A^{(i)}) \wdt\th_k + o(\wdt\th_k) \r^2 I_{\{\al_k = a_i\}}
\\
&\ \le E \sum^\mz_{i=1}E^\mu_t K \mu \sum^{(t+s)/\mu -1}_{k=t/\mu}
\sum^{(t+s)/\mu-1}_{j=t/\mu}
[(A_k^{(i)} - A^{(i)}) \wdt\th_k]'[(A_j^{(i)} - A^{(i)})\wdt\th_j ]
I_{\{\al_k = a_i=\al_j\}}
\\
& \quad + K\mu \sum^{(t+s)/\mu -1}_{k=t/\mu} E \l \wdt\th_k \r^2
\\
&\ \le E \sum^\mz_{i=1}K\mu \Big[  E^\mu_t
\sum^{(t+s)/\mu -1}_{k=t/\mu} \l (A_k^{(i)} - A^{(i)}) \sqrt\mu u_k \r
\sum_{j\ge k}\l (A_j^{(i)} - A^{(i)}) \sqrt\mu u_j \r \Big] I_{\{\al_k =
a_i=\al_j\}} +O(\mu)s
\\
&\ \le O(\mu)s
}
For any $T<\infty$ and any $0\le t\le T$,
$$\lim_{\delta \to 0} \limsup_{\mu \to 0} \Big\{\sup_{0\le s\le \delta}
 E [E^\mu_t |u^\mu(s+t) - u^\mu(t)|^2]\Big\} =
0,$$ so $\{u^\mu \cd\}$ is tight.
\qed

Note that $g_k(a_i, i) = \ph_k \sg(\ph'_k[a_i - a_i] +e_k) = \ph_k
\sg(e_k)$.
The following is a variant of the well-known central limit theorem for
mixing processes;
see \cite{Billingsley} or \cite{Ethier} for details.

\begin{lem}\label{bromo}
Define $\varpi_k := \ph_k \sg(e_k)$. Then
\beq{BM-1}
\sqrt \mu \sum^{(t/\mu)-1}_{j=0} \varpi_j \ \hbox{ converges weakly to a
Brownian motion } \ \wdt w(t)
\eeq
with covariance $\wdt \Sigma t$
such that the covariance $\wdt \Sigma$ is given by
\beq{BM-2}
\wdt \Sigma = E \varpi_0 \varpi'_0 + \sum^\infty_{j=1}E \varpi_j
\varpi_0'
+ \sum^\infty_{j=1} E \varpi_0 \varpi_j'.
\eeq
\end{lem}

\begin{thm}\label{rate-e}
$u^\mu\cd$ converges weakly to $u\cd$ such that
$u\cd$ is the solution of
\beq{sde} d u(t)= -A^{(\al)} u dt - \wdt \Sigma^{1/2} dw,\eeq
where $w\cd$ is a standard Brownian motion.
\end{thm}

\para{Proof.}
As usual, extract a convergent subsequence of $u^\mu \cd$ (still denoted
by $u^\mu \cd$)
with limit $u \cd$. We will show that for each $s,t >0$, the limit
process satisfies
\beq{sde-2}
u(t+s) - u(t) = \int^{t+s}_t -A^{(\al(\tau)}u(\tau)d\tau
- \int^{t+s}_t \wdt \Sigma^{1/2} dw(\tau)
\eeq
Note from (\ref{udiff-1}),
\aleq{\label{udiff-mu}
\disp u^\mu(t+s) - u^\mu(t) & =
-\sqrt \mu \sum^{(t+s)/\mu -1}_{k=t/\mu}g_k
+ O(\sqrt \mu)
\\
& = \sum^\mz_{i=1} [-\sqrt \mu \sum^{(t+s)/\mu -1}_{k=t/\mu}g_k]
I_{\{\al_k= a_i\}}
+ O(\sqrt\mu)
.}
Define \bea \ad g_k(i) := g_k I_{\{\al_k= a_i\}}, \ \wdt g_k(i) := E_k g_k(i),
\ \hbox{ and } \\
\ad \Delta_k(i) := [g_k(i) -g_k(a_i, i) -(\wdt g_k(i) -\wdt g_k(a_i, i))].\eea
We then expand on the (negative of the) inside of the sum indexed by
$i$ in (\ref{udiff-mu}) as
\aleq{\label{g-expd}
 &\!\!\!\sqrt \mu \sum^{(t+s)/\mu -1}_{k=t/\mu} g_k(i)
\\
&= \sum^{(t+s)/\mu -1}_{k=t/\mu}\sqrt\mu g_k(a_i, i) +
\sum^{(t+s)/\mu -1}_{k=t/\mu}\sqrt\mu[\wdt g_k(i)-\wdt g_k (a_i, i)]
 + \sum^{(t+s)/\mu -1}_{k=t/\mu}\sqrt\mu \Delta_k(i)
\\
& = \sum^{(t+s)/\mu -1}_{k=t/\mu} \sqrt\mu \varpi_k
 + \sum^{(t+s)/\mu -1}_{k=t/\mu} \mu [A^{(i)}_k u_k +o(|u_k|)]
 + \sum^{(t+s)/\mu -1}_{k=t/\mu}\sqrt\mu \Delta_k(i)
.}
Note that for the second term above we used
$\wdt g_k(a_i, i) = o(\wdt\th_k) = o(\sqrt\mu |u_k|)$ by {\rm(A3)}.
First, we show the last term in (\ref{g-expd}) is $o(1)$. Since
$\Delta_k(i)$ is a martingale difference, we have
\aleq{\label{Delta-es}
\disp E \l \sum^{(t+s)/\mu -1}_{k=t/\mu}\sqrt\mu \Delta_k(i) \r^2
& =  \sum^{(t+s)/\mu -1}_{k=t/\mu}\mu \l \Delta_k(i) \r^2
\\
& =  \sum^{(t+s)/\mu -1}_{k=t/\mu}\mu E [g_k(i) - g_k(a_i, i)]'[g_k(i)
- g_k(a_i, i)]
\\
& \qquad\qquad + \sum^{(t+s)/\mu -1}_{k=t/\mu}\mu E [\wdt g_k(i) - \wdt g_k(a_i,
i)]'[\wdt g_k(i) - \wdt g_k(a_i, i)]
.}
The boundedness of $\ph_k$ and $u_k$ implies $\sqrt\mu \ph'_k u_k \to
0$ in probability 
uniformly in $k$
as $\mu \to 0$. Hence, the first term in (\ref{Delta-es}) has
\aleq{\label{Delta-es1}
 & \sum^{(t+s)/\mu -1}_{k=t/\mu}\mu E [g_k(i) - g_k(a_i, i)]'[g_k(i) -
g_k(a_i, i)]
\\& = \sum^{(t+s)/\mu -1}_{k=t/\mu}\mu E\ph'_k\ph_k [\sg(\sqrt\mu \ph'_k
u_k +e_k) -\sg(e_k)]^2
 \to 0 \ \hbox{ as } \mu \to 0.
}
Using {\rm (A3)} and {\rm (A4)}, along with the boundedness of $u_k$, on
the second term of (\ref{Delta-es}) gives
\aleq{\label{Delta-es2}
&\!\!\! \sum^{(t+s)/\mu -1}_{k=t/\mu}\mu E [\wdt g_k(i) - \wdt g_k(a_i,
i)]'[\wdt g_k(i) - \wdt g_k(a_i, i)]
\\
&\ = \sum^{(t+s)/\mu -1}_{k=t/\mu}\mu E[\sqrt\mu A^{(i)}_k u_k
+ o(\sqrt\mu |u_k|)]' [\sqrt\mu A^{(i)}_k u_k + o(\sqrt\mu |u_k|)]
\\
&\ = \sum^{(t+s)/\mu -1}_{k=t/\mu} \mu^2
E \l(A^{(i)}_k - A^{(i)})u_k + A^{(i)}u_k + o(|u_k|) \r^2
\\
&\ \le \mu^2K \sum^\infty_{k=t/\mu} \phi(k-t/\mu)
+ \mu^2 \sum^{(t+s)/\mu -1}_{k=t/\mu} K
 \to 0 \ \hbox{ as } \mu \to 0.
}
Hence $$E  \sum^{(t+s)/\mu -1}_{k=t/\mu} \sqrt\mu \Delta_k(i) \to 0 \
\hbox{ as } \mu \to 0.$$

Next, in the second term of (\ref{g-expd}) we have
\aleq{\label{drift-1}
&\!\! \sum^{(t+s)/\mu -1}_{k=t/\mu} \mu [A^{(i)}_k u_k +o(|u_k|)]\\
& \quad = A^{(i)} \sum^{(t+s)/\mu -1}_{k=t/\mu} \mu u_k
+\sum^{(t+s)/\mu -1}_{k=t/\mu} \mu (A^{(i)}_k - A^{(i)})u_k
+\sum^{(t+s)/\mu -1}_{k=t/\mu} \mu o(|u_k|)
.}
Similar to the previous section, choose a sequence $m_\mu$ such that
$m_\mu \to \infty$
as $\mu \to 0$ but $\dl_\mu / \sqrt\mu = \sqrt\mu m_\mu \to 0$. Then
\aleq{\label{drift-2}
\sum^{(t+s)/\mu -1}_{k=t/\mu} \mu u_k
= \sum^{t+s}_{l \dl_\mu=t} \dl_\mu u_{l m_\mu}
+ \sum^{t+s}_{l \dl_\mu=t}\dl_\mu {1 \over m_\mu}
\sum^{lm_\mu+m_\mu-1}_{k=lm_\mu}[u_k-u_{l m_\mu}]
.}
Since for $l m_\mu \le k < lm_\mu + m_\mu$,
$u_k - u_{lm_\mu} = O(\dl_\mu / \sqrt\mu)$,
so the second term above goes to 0 in probability, uniformly in t.
Similarly, by {\rm(A3)},
\aleq{\label{drift-3}
& \sum^{(t+s)/\mu -1}_{k=t/\mu} \mu (A^{(i)}_k - A^{(i)})u_k
= \sum^{t+s}_{l \dl_\mu=t}\dl_\mu {1 \over m_\mu}
\sum^{lm_\mu+m_\mu-1}_{k=lm_\mu} (A^{(i)}_k - A^{(i)})u_k
\to 0
.}
Likewise, $\sum^{(t+s)/\mu -1}_{k=t/\mu} \mu o(|u_k|) \to 0$ in
probability uniformly in t.

Hence, putting the above estimates together we obtain
\aleq{\label{sde-3}
\disp u(t+s) - u(t) & = \lim_{\mu \to 0} u^\mu(t+s) - u^\mu(t)
 \\
 & = \lim_{\mu \to 0} \sum^\mz_{i=1} \Big[
-A^{i} \sum^{t+s}_{l \dl_\mu=t} \dl_\mu u_{l m_\mu} \Big] I_{\al_k = a_i}
- \sum^{(t+s)/\mu -1}_{k=t/\mu} \sqrt\mu \varpi_k
\\
&= -\int^{t+s}_t A^{\al(\tau)}u(\tau)d\tau -\int^{t+s}_t
\wdt\Sigma^{1/2}dw(\tau)
.}
\qed

\subsection{Scaled Errors: $\e \ll \mu$}
The analysis for the cases $\e \ll \mu$ and $\e \gg \mu$ is similar to that for $\e=O(\mu)$.
We omit the details and present the main results. Recall that in the $\e \ll \mu$ case, the parameter is essentially a constant and thus we look to the initial distribution to determine the asymptotic properties.

Define
\begin{align*}
 &\al_{*} := \sum_{i=1}^\mz a_i P(\al_0 = a_i),
\qquad v_n := {\al_* - \th_n \over \sqrt \mu }
\\& v^\mu(t):=v_n \quad  \hbox{for} \quad t\in [ (n-N_\mu)\mu, (n-N_\mu)\mu+\mu)
\\& A^{(*)} := \sum_{i=1}^\mz A^{(i)} P(\al_0 = a_i)
\end{align*}
Then we have the following:

\begin{thm}\label{rate-sl}
Assume $\e=\mu^{1+\Delta}$ for some $\Delta > 0$. Then
$v^\mu\cd$ converges weakly to $v\cd$ such that
$v\cd$ is the solution of
\beq{sde} d v(t)= -A^{(*)} v dt - \wdt \Sigma^{1/2} dw,\eeq
where $w\cd$ is a standard Brownian motion.
\end{thm}

\subsection{Scaled Errors: $\e \gg \mu$}
Again, here the idea is that the parameter varies so quickly that it quickly converges to the stationary distribution $\nu=(\nu_1,\dots,\nu_{m_0})$. Thus we look to the expectation against the stationary distribution to determine the asymptotic properties.

Define
\begin{align*}
 &\bar{\al} := \sum_{i=1}^\mz a_i \nu_i
\qquad z_n := {\bar{\al} - \th_n \over \sqrt \mu }
\\& z^\mu(t) := z_n \quad  \hbox{for} \quad t\in [ (n-N_\mu)\mu, (n-N_\mu)\mu+\mu)
\\& \bar{A} := \sum_{i=1}^\mz A^{(i)} \nu_i
\end{align*}
We have the following result.

\begin{thm}\label{rate-fa}
Assume $\e=\mu^{\gamma}$ for some $1/2 < \gamma <1$. Then
$z^\mu\cd$ converges weakly to $z\cd$ such that
$z\cd$ is the solution of
\beq{sde} d z(t)= -\bar{A} z dt - \wdt \Sigma^{1/2} dw,\eeq
where $w\cd$ is a standard Brownian motion.
\end{thm}

\section{Numerical Examples}\label{sec:num}

Here we demonstrate the performance of the Sign-Error (SE) algorithms and compare it with the Sign-Regressor (SR) and Least Mean Squares (LMS) algorithms (see \cite{YHW11, YinK05} respectively). We fix the step size $\mu = .05$ and consider three cases: $\e = (3/5) \mu$ ( $\e = O(\mu)$); $\e = \mu^2$ (a slowly-varying Markov chain); and $\e=\sqrt\mu$ (a fast Markov chain).

We use state space $\mathcal{M} =\{-1,0, 1\}$ with transition matrix $P^\e = I + \e Q$, where
$$Q =
\left[ {\begin{array}{rrr}
-0.6 & 0.4 & 0.2 \\
0.2 & -0.5 & 0.3 \\
0.4 & 0.1 & -0.5 \\
\end{array} } \right].$$
is the generator of a continuous-time Markov chain whose stationary distribution is therefore
$\nu = (1/3, 1/3, 1/3) $. Hence $ \lbar\al = \sum_{i=1}^3 a_i \nu_i= 0$
We take the initial distribution for $\al_0$ to be $(3/4,1/8, 1/8)$. So
$\al_* = \sum_{i=1}^3 a_i P( \al_0 =a_i) = -0.625$.
$\{\ph_n \}$ and $\{e_n\}$ are  i.i.d.
 $\mathcal{N} (0,1)$ and $\mathcal{N} (0,.25)$, respectively. We proceed to observe $1000$ iterations of the algorithm for the cases $\e = O(\mu)$ and $\e >> \mu$, and $10,000$ iterations for the case $\e << \mu$ (in order to illustrate some variations of the parameter).

\begin{figure}[htb]
\centerline{\epsfxsize=5in \epsfysize=3in \epsfbox{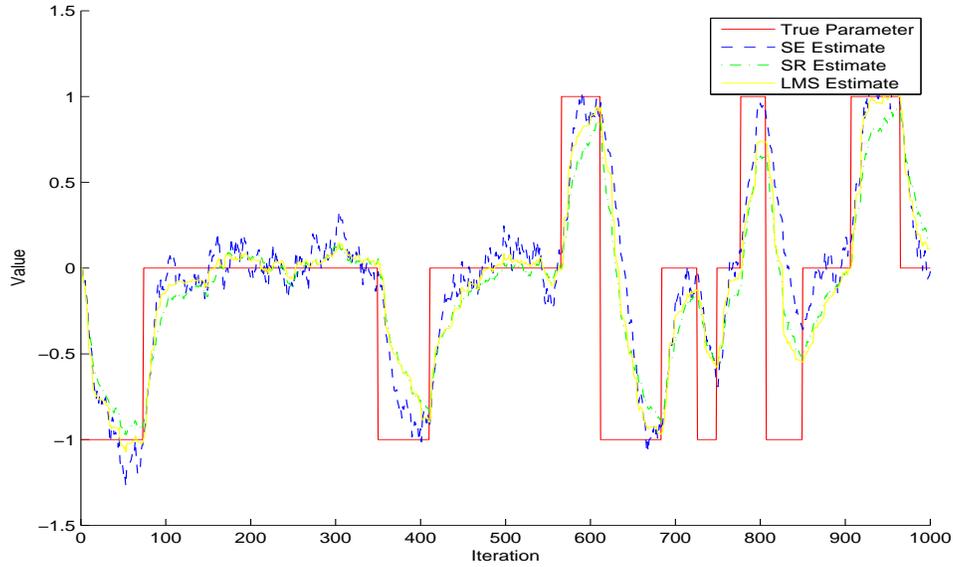}}
\caption{Markov chain parameter process and
 estimates obtained by adaptive filtering with  $\e = O(\mu)$}
 \label{e=mu}
\end{figure}

\begin{figure}[h!tb]
\begin{minipage}[b]{0.5\linewidth}
\centering
\includegraphics[scale=.3]{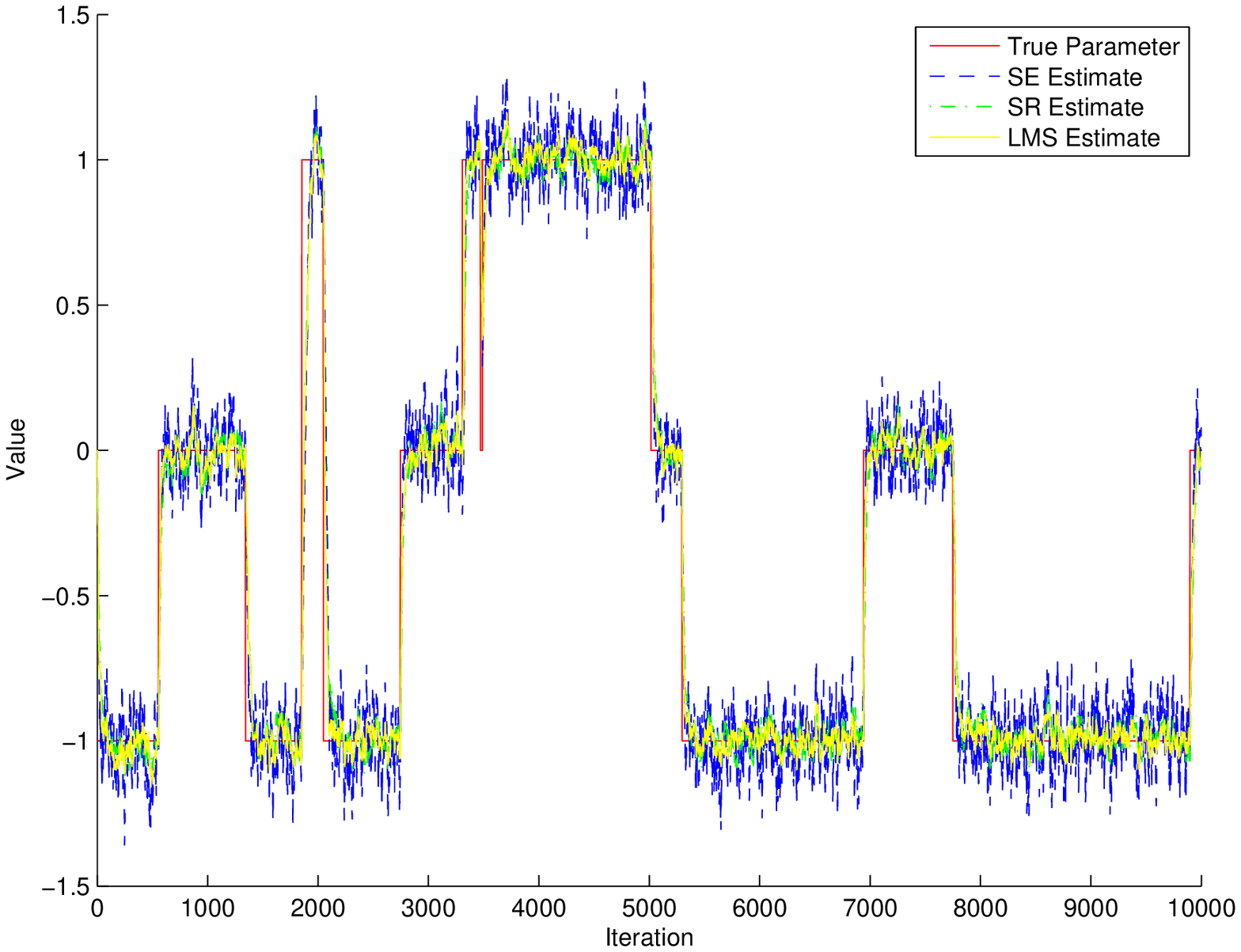}
\caption{Markov chain parameter process and
 estimates obtained by adaptive filtering with  $\e \ll O(\mu)$}
 \label{elmu}
\end{minipage}
\hspace{0.5cm}
\begin{minipage}[b]{0.5\linewidth}
\centering
\includegraphics[scale=.3]{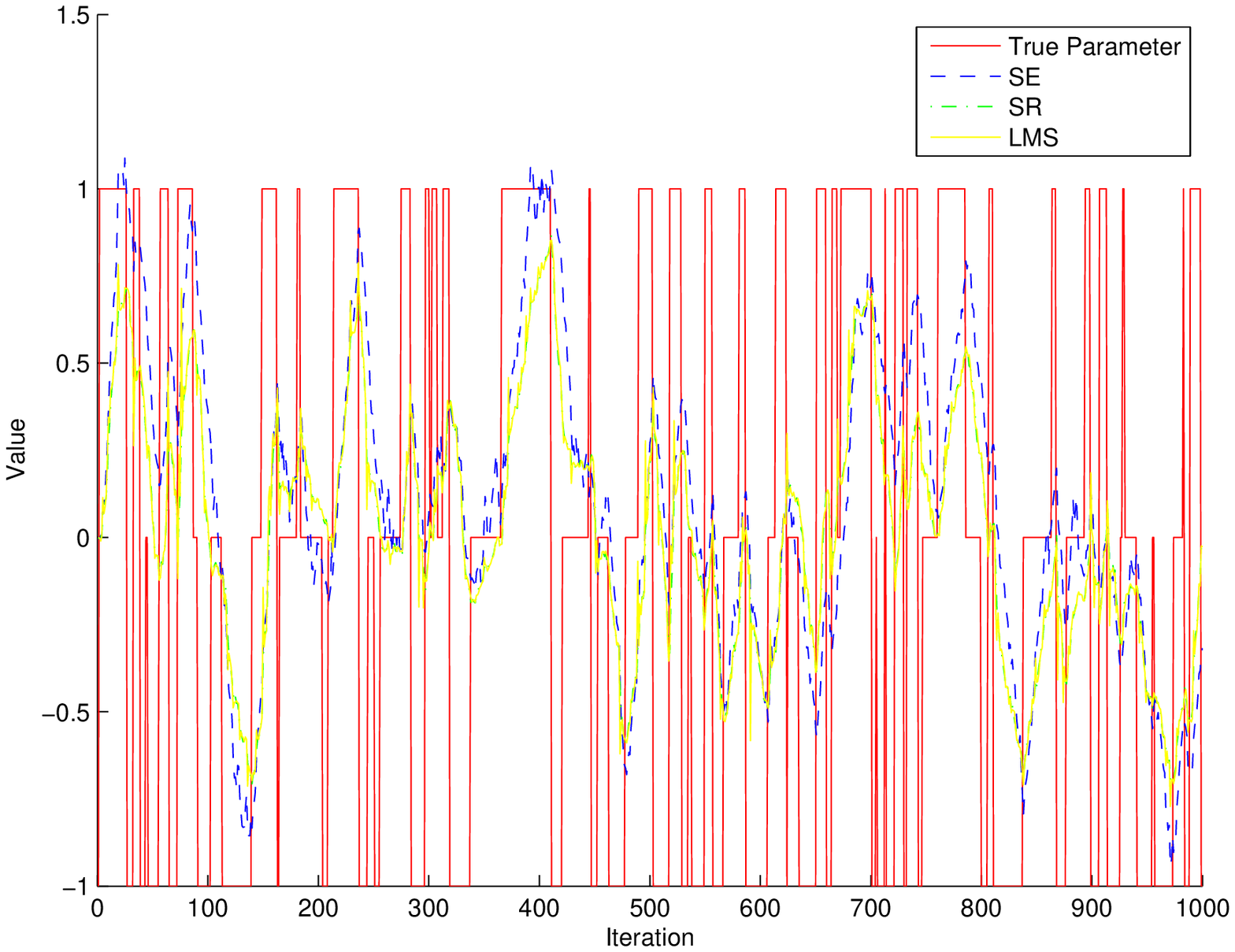}
\caption{Markov chain parameter process and
 estimates obtained by adaptive filtering with  $\e \gg \mu$ }
 \label{egmu}
\end{minipage}
\hspace{0.5cm}
\end{figure}

\begin{figure}[h!tb]
\begin{minipage}[b]{0.5\linewidth}
\centering
\includegraphics[scale=.3]{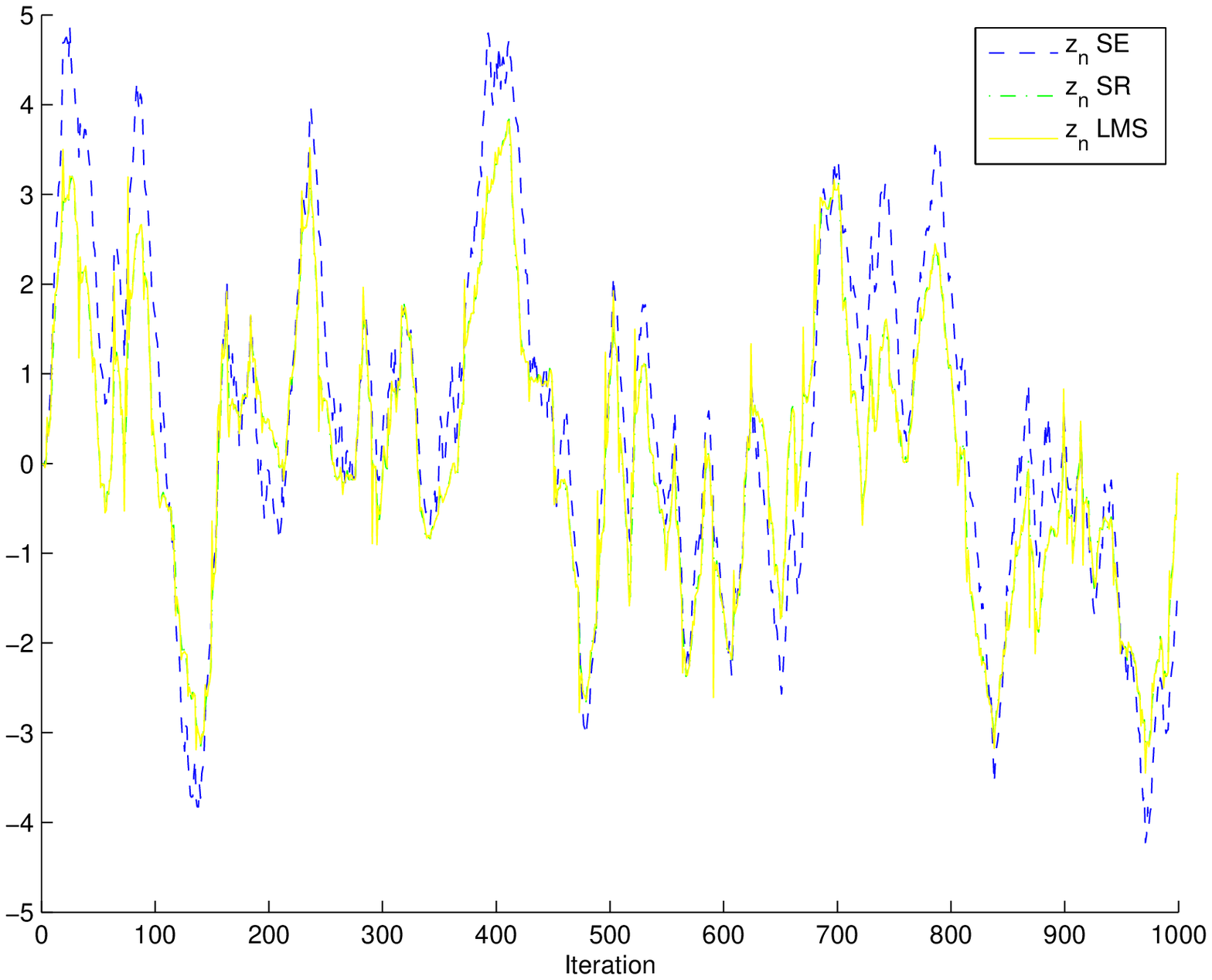}
\caption{Scaled error $z_n$ with
fast varying Markov chain ($\e \gg \mu$):
Diffusion behavior}
 \label{diff}
\end{minipage}
\hspace{0.5cm}
\begin{minipage}[b]{0.5\linewidth}
\centering
\includegraphics[scale=.3]{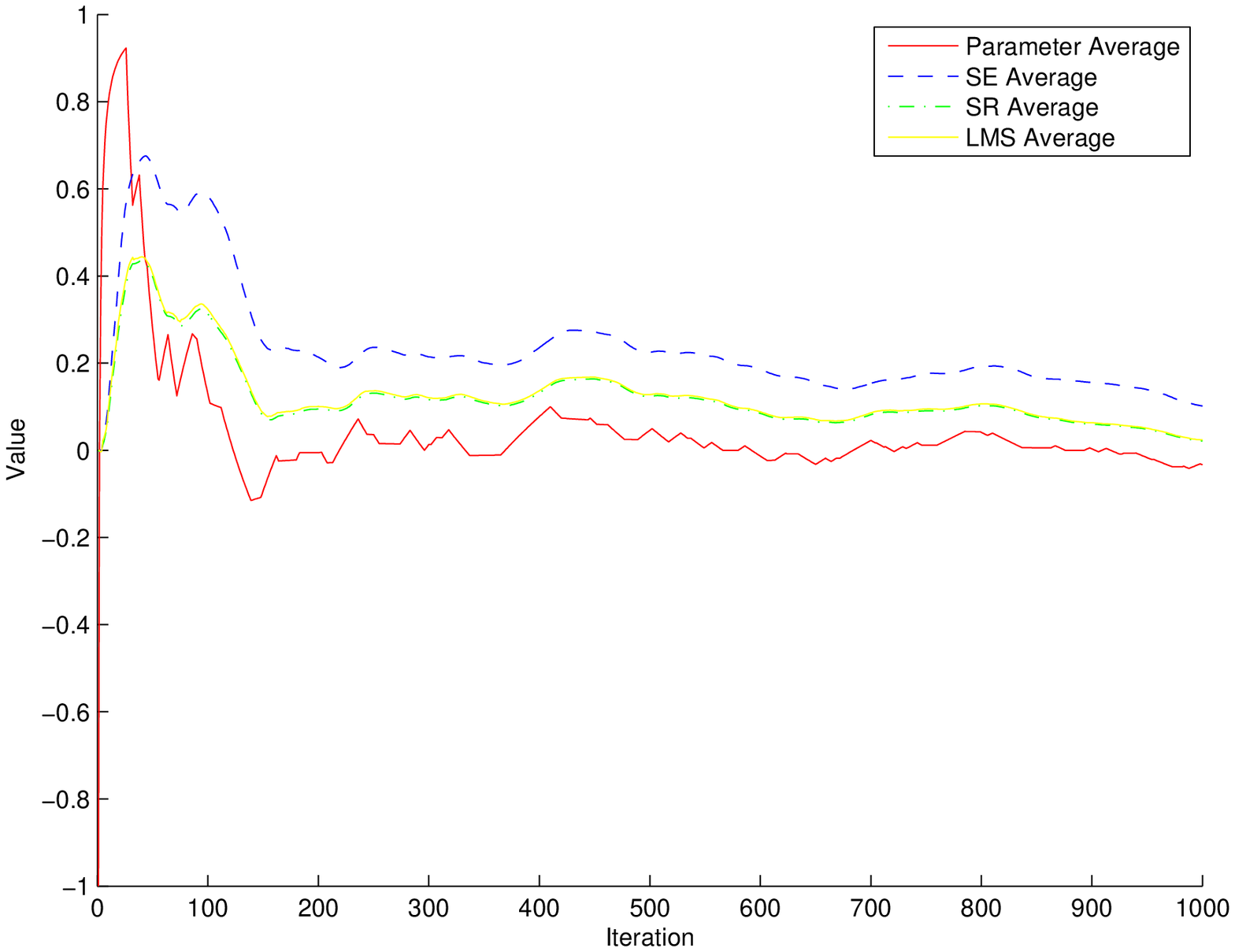}
\caption{Average of parameter process and estimates over time with $\e \gg \mu$}
 \label{egmu}
\end{minipage}
\hspace{0.5cm}
\end{figure}

To observe the tracking behavior of the SE algorithm, in comparison to the SR and LMS algorithms, we overlay the respective plots for each case.
When $\e = O(\mu)$, the LMS and SR estimates tend to be approximately equal, while the SE estimates show more deviations from the other estimates. The SE algorithm responds to changes in the parameter more quickly, while the LMS and SR algorithms adhere to the parameter more closely while it is stationary. In the $\e \ll \mu$ case, we see this behavior repeated. While all three estimates track the parameter closely, the LMS and SR estimates deviate from the parameter less than the SE estimates between jumps of the parameter.

 In the $\e \gg \mu$ case, none of the algorithms can track the parameter at each iterate very well. However, when we observe the scaled error against the stationary distribution of the Markov chain $z_n$,
the expected diffusion behavior is displayed. Examining the cumulative average of the parameter and the estimates of the iterates, we note that the  parameter average quickly converges to
$\bar{\al}$. The LMS and SR estimate averages adhere closely to the parameter average, while the SE estimate average deviates slightly more.

\FloatBarrier

\end{document}